\documentclass[11pt]{article}
\usepackage{amsmath,amsthm,amsfonts,amscd,bezier}
\usepackage{amssymb}
\usepackage[usenames]{color}
\usepackage[T1]{fontenc}       
\usepackage[all]{xy}

\newtheorem{lemma}{Lemma}[section]
\newtheorem{example}[lemma]{Example}
\newtheorem{definition}[lemma]{Definition}
\newtheorem{corollary}[lemma]{Corollary}
\newtheorem{remark}[lemma]{Remark}
\newtheorem{theorem}[lemma]{Theorem}
\newtheorem{proposition}[lemma]{Proposition}
\newcommand{\Dem}{\noindent{\sc Proof:\ \ }}
\newcommand{\cqd}{{\hfill $\rule{1.5mm}{1.5mm}$}\vspace{0.2cm}}

\def\Cp{\mathbb{C}}
\def\N{\mathbb{N}}
\def\Cm{\rm{\mathbb{C}}}
\def\O{{\cal O}}
\def\C{{\cal C}}

\def\M{{\cal M}}

\def\Tg{T{\cal A}_e(f)}
\def\A{{\mathcal A}}
\def\L{{\mathcal L}}
\def\Acod{{\cal A}_ecod(f)}

\begin{document}

\pagestyle{plain} \pagenumbering{arabic}

\title{Join operation and $\A$-finite map-germs}
\author{Rodrigues Hernandes, M. E.  and Ruas, M. A. S. \thanks{The first author was partially supported by FAPESP Proc. $2019/07316-0$ and the second author by FAPESP Proc. $2019/21181-0 $ and CNPq Proc. $305695/2019-3$.}}
\date{} \maketitle

\begin{abstract}
In this work we define some map-germs, called elementary joins, for the purpose of producing new $\A$-finite map-germs from them. In particular, we describe a general form of an $\A$-finite monomial map from $(\Cp^n,0)$ to $(\Cp^{p},0)$ for $p\geq 2n$ of any corank in terms of elementary join maps. Our main tools are the delta invariant and some invariants of curves.
\end{abstract}

\vspace{0.3cm}

\begin{center}
\end{center}

\begin{center}
 \end{center}

\vspace{0.5cm}

\section{Introduction}

The class of map-germs from $(\mathbb{K}^n,0)$ to $(\mathbb{K}^{2n},0)$, where $\mathbb{K}$ is $\mathbb{R}$ or $\Cp$, has been investigated by many authors using different approaches. For instance, Gaffney gives in \cite{Gaffney} a formula for the number of double points at the origin of
a stable perturbation of a map-germ $f:(\Cp^n,0) \to (\Cp^{2n},0)$
in terms of the zero dimensional Segre number of the double point
ideal of $f$ and the number of Whitney umbrella points of the
composition of $f$ with a generic projection from $\Cp^{2n}$ to
$\Cp^{2n-1}$. Klotz, Pop and Rieger classify $\A$-simple germs $(\mathbb{R}^n,0)\to
(\mathbb{R}^{2n},0)$ (see \cite{Klotz}) to show that all of them have a
real deformation with the maximal number of real double-points. Symplectic invariants of curves and surfaces in this class of maps were also studied in \cite{Ishikawa} by Ishikawa and Janeczko. Parameterized singular surfaces in $\mathbb{R}^4$ appear in \cite{BMN}, where Birbrair, Mendes and Nu\~{n}o-Ballesteros study relations between topological and metric properties, and also in \cite{BRS} where Benedini, Ruas and Sinha study their second order geometry. In \cite{RodRuas}, the authors classify $\A$-finitely determined monomial surfaces in $\Cp^4$.

Our main purpose is to present an ``operation'', called join operation, among a finite number of map-germs in order to produce classes of $\A$-finitely determined (or $\A$-finite) map-germs from $(\Cp^n,0)$ to $(\Cp^p,0)$ with $p \geq 2n$. In particular, we provide families of these maps of any corank. For two map-germs $F$ and $G$, this operation consists in joining them by means of a germ $H$, denoted $(F*G)_H$, defined as $(F*G)_H(x,y):=(F(x),G(y),H(x,y))$ which preserve in some sense the properties of $F$ and $G$.  Different type of operations to obtain new $\A$-finite germs from $\Cp^n$ to $\Cp^p$ from germs in lower dimensions were defined in \cite{CMW}, see also \cite{ORW}.

We define basic map-germs, called elementary joins, which are obtained by joining the identity and a parameterized curve or a finite number of parameterized curves. These maps are $\A$-finite (see Propositions \ref{FamilyA-Finite} and \ref{PropCorankn-map}). In Theorem \ref{ClassMonomProjCurve} and Corollaries \ref{ClassifMon-n2n} and \ref{ClassMon-General} we obtain normal forms of some classes of $\A$-finite monomial map-germs which are completely characterized in terms of the elementary joins. The result in Corollary \ref{ClassifMon-n2n} generalizes for map-germs from $\Cp^n$ to $\Cp^{2n}$ the classification given in Theorem 2.4 in \cite{RodRuas} for surfaces in $\Cp^4$. We also give the general form of an $\A$-finite monomial map-germ $f:(\Cp^n,0) \to (\Cp^{p},0)$ for $p\geq 2n$ of any corank. More precisely, in Theorem \ref{ClassifMonoCurve-Corankq}, we prove that if $f$ is a corank $q$ monomial map-germ, then $f$ is $\A$-finite if and only if $f$ is $\A$-equivalent to the join of the identity in $\Cp^{n-q}$ with $q$ parameterized curves. As a corollary we obtain infinitely many families of $\A$-finite map-germs, not necessarily monomial. Two fundamental tools to obtain the previous results were the invariants: the semigroup associated to a parameterized curve and the delta invariant $\delta_f:=dim_{\Cp}\frac{\O_n}{f^*(\O_p)}$, defined in \cite{RodRuas}.

The monomial normal forms of Theorem \ref{ClassifMonoCurve-Corankq} can be seen as $\A$-finite reflection maps, as defined in \cite{P}. In addition, their images $f(\Cp^n)$ are non normal toric varieties in $\Cp^p$ with isolated singularity at the origin, as they are the loci of a binomial ideal (see \cite{ES}). Conditions for the equisingularity of families of isolated non-normal singularities parameterized by this class of map-germs were investigated in \cite{Greuel}.

In Section 5, we focus our study on some $\A$-invariants of map-germs $f:(\Cp^n,0) \to (\Cp^{p},0)$ for $p\geq 2n$. The delta invariant is also useful for verifying whether a map is $\A$-finitely determined, since $\delta_f$ is easier to be calculated than the $\A_e$-codimension of $f$ ($\A_ecod(f)$). If $p=2n$, we also have the number $d(f)$ of double points appearing in a stable perturbation of $f$. The results of this section provide estimates of $\delta_f$ and consequently for $\Acod$ and some formulas for $d(f)$.  In Subsection 5.1, we consider $f$ $\A$-equivalent to $(g_k,h)$, where $g_k=\pi_k\circ f$ in which $\pi_k:\Cp^p \to \Cp^k$ is a canonical projection and $h:(\Cp^n,0) \to (\Cp^{p-k},0)$. Such decomposition can be interesting in the study of analytic invariants, if one knows some properties of $g_k$ and $h$.  In Theorem \ref{AcodDelta}, we prove that $\A_ecod(f)\leq \A_ecod(g_k)+(p-k)\delta_f$, with a formula when $g_k$ is a stable map.  If the projection $g_n$ of $f$ is a fold map, we prove in Proposition \ref{DobraDelta} that $d(f)$ and $\delta_f$ coincide and we provide upper and lower bounds for $\A_ecod(f)$.

\section{Basic notations of singularities}
\label{Sec-Notations}
Analytical invariants play a central role in this work. The delta invariant of a map-germ and some invariants of the theory of curves such as the semigroup of values and its set of gaps are fundamental for the characterization of some classes of map-germs.

Let $\O(n,p)$ be the set of smooth map-germs $f:(\Cp^n,0) \to
\Cp^p$, which is a free $\O_n$-module of rank $p$, where $\O_n$ is
the local ring of function germs in $\Cp^n$ at the origin, with
maximal ideal $\M_n$. Denoting by ${\mathcal R}=Diff(\Cp^n,0)$ and ${\mathcal L}=Diff(\Cp^p,0)$ the groups of germs of diffeomorphisms at the origin, we consider the set $\M_n\O(n,p)=\{f\in \O(n,p);f(0)=0\}$ under the action of the group $\A={\mathcal R}\times {\mathcal L}$. Two map-germs $f,g\in \M_n\O(n,p)$ are called $\A$-equivalent, denoted $f \underset{\A}{\sim} g$, if there exist $(\rho,\psi)\in \A$ such that $g=\psi \circ f \circ \rho^{-1}$.

The set of vector fields along $f$ is defined by $\theta(f)=\{\sigma: (\Cp^n,0) \to T\Cp^p; \
\pi \circ \sigma= f\}$ such that $\pi: T\Cp^p\to (\Cp^p,0)$ is the natural projection. It is a free $\O_n$-module of rank $p$. The set of germs of vector fields in $\Cp^n$ at the origin is a free $\O_n$-module of rank $n$, denoted by $\theta(n)$. Similarly, $\theta(p)$ denotes the $\O_p$-module of rank $p$ of vector fields in $\Cp^p$ at the origin. We define  $tf: \theta(n) \to
\theta(f)$ by $tf(\epsilon)(x)=df_x(\epsilon(x))$, where $df$ is the differential of $f$ and $wf: \theta{(p)}
\to \theta(f)$ such that $wf(\eta) =f^*(\eta):=\eta \circ f$. The {\it extended tangent space} $T{\A}_e(f)$ to the $\A$-orbit $\A(f)$ at the point $f$ is the complex vector space $T{\A}_e(f)=tf(\theta(n)) + wf(\theta{(p)})$.  The \textit{extended tangent space} $T{\mathcal L}_e(f)$ to the ${\mathcal L}$-orbit of $f$ is defined as $T{\mathcal L}_e(f)= wf(\theta{(p)})$.

If ${\mathcal G}$ is one of the groups of Mather ${\mathcal K}, \A,{\mathcal L}$ or ${\mathcal R}$, a smooth map-germ $f:(\Cp^n,0) \to
(\Cp^p,0)$ is called ${\mathcal G}$-finitely determined (or ${\mathcal G}$-finite) if there exists $k < \infty$ such that for all smooth map-germs $g:(\Cp^n,0) \to
(\Cp^p,0)$ with $j^kf(0)=j^kg(0)$ we have that $f$ is ${\mathcal G}$-equivalent to $g$, where $j^kf(0)$ denotes the $k$-jet of $f$ at the origin.

 Finitely determined map-germs in $\M_n\O(n,p)$ are particularly interesting for the class $p\geq 2n$. In this class of map-germs, the Mather-Gaffney geometric criterion of finite determinacy states that $f$ is $\A$-finite if and only if there exists a neighborhood $U$ at the origin in $\Cp^n$ such that $0$ is an isolated singularity of $f$ and $f|_{U\setminus \{0\}}$ is one-to-one. Moreover, Gaffney proves that $f$ is $\A$-finite if and only if $f$ is ${\cal L}$-finite (see \cite[Theorem 2.5]{Wall}). This result can be algebraically translated in terms of the ${\cal G}_e$-codimension of $f$, defined as ${\cal G}_ecod(f)=dim_{\Cp}\frac{\theta(f)}{T{\cal G}_e(f)}$, where ${\cal G}=\A$ or ${\cal L}$. More precisely, by the Infinitesimal criterion of finite determinacy (\cite[Theorem 1.2]{Wall}), $f$ is ${\mathcal G}$-finite if and only if ${\cal G}_e$-codimension of $f$ is finite. Therefore, if $p\geq 2n$ then $\Acod$ is finite if and only if ${\cal L}_ecod(f)$ is finite. 

In \cite{RodRuas} the authors consider the so-called delta invariant of $f$ defined by $\delta_f:=dim_{\Cm}
\frac{\O_n}{f^*(\O_{p})}$. Notice that
$$\frac{\theta(f)}{wf(\theta{(p)})}\ \simeq \ \frac{\oplus_{i=1}^p \O_n}{f^*\left(\oplus_{i=1}^p \O_p\right)} \ \simeq \ \bigoplus_{i=1}^p \frac{\O_n}{f^*(\O_{p})},$$ then it follows that ${\cal L}_ecod(f)=dim_{\Cp}\ \frac{\theta(f)}{wf(\theta{(p)})}=p \cdot \delta_f$. Thus, for $p\geq 2n$ we obtain that $f$ is $\A$-finite if and only if $\delta_f$ is finite. In other words, the finite determinacy of $f$ is described in terms of the delta invariant, which is easier to be calculated than $\Acod$.  By a previous relation $\delta_f$ is an $\L$-invariant. More strongly, it is an $\A$-invariant. In fact, Gaffney in his thesis (see \cite{Wall}, Theorem 2.7(iv)) proves that if $f,g \in \O(n,p)$ are $\A$-finite and $n < p$, then $f \underset{\A}{\sim} g$ if and only if the algebras $f^*(\O_p)$ and $g^*(\O_p)$ are isomorphic.

\vspace{0.2cm}
In what follows we introduce some tools of the theory of curves. For more details see \cite{Hefez}.

Let $\phi_k:(\Cp,0)\to (\Cp^k,0)$ be a parameterized curve in $\Cp^k$ in which we assume (up to $\A$-equivalence) be given in the form
\begin{equation}
\label{ParamCurve}
\phi_k(y)=(y^{m_1}, y^{m_2}+\sum_{i>m_2}a_{i2}y^{i},\ldots, y^{m_k}+\sum_{i>m_k}a_{ik}y^{i})
\end{equation}
with $1 < m_1\leq m_2 \leq \ldots \leq m_k$. We suppose that the map $\phi_k$ is primitive, that is, we cannot reparameterize it by a power of a new variable, or equivalently the greatest common divisor of the exponents of all components of $\phi_k$ is equal to $1$. The integer $m_1$ is called the \textit{multiplicity} of the curve.

\begin{itemize}
	\item {\sc Numerical Semigroups and Ap\'{e}ry sequence}
	
	A subset $\mathbf{S}\neq \{0\}$ is called a semigroup of $\N$ if $\mathbf{S}$ contains the element $0$ and it is closed under addition. Any semigroup $\mathbf{S}$ of $\N$ is finitely generated, that is, there exist $v_1,\ldots, v_r \in \mathbf{S}$ such that $\mathbf{S}:=\langle v_1, \ldots, v_r \rangle=\{\lambda_1v_1+\ldots + \lambda_r v_r;\ \lambda_i \in \N, i=1, \ldots, r\}.$
	If $r \geq 2$ and $gcd(v_1,\ldots, v_r)=1$, then $\mathbf{S}$ has a conductor, which means that there exists an element $c\in \mathbf{S}$ such that $c-1 \notin \mathbf{S}$ and $c+\N \subset \mathbf{S}$. Consequently  $\N \setminus \mathbf{S}$ is finite and its elements are called {\it gaps} of $\mathbf{S}$.
	
 In order to describe the gaps of $\mathbf{S}$ we introduce the notion of the {\it Ap\'{e}ry sequence}  of $\mathbf{S}$ with respect to an element $q \in \mathbf{S}\setminus \{0\}$, defined inductively by: $a_0=0$ and for $j=1, \ldots, q-1$,
\begin{equation}
\label{AperySet}
a_j=min\left(\mathbf{S}\setminus \bigcup_{i=0}^{j-1}(a_i+q\N)\right).
\end{equation}

For more details about these concepts see \cite{Hefez}.

\item {\sc Semigroup and gaps associated to a curve}

Let $(\C,0)\subset (\Cp^k,0)$ be an analytic irreducible curve, parameterized by $\phi_k$ as in (\ref{ParamCurve}), with $\O:=\Cp\{X_1,\ldots, X_k\}/I$ the local ring of $(\C,0)$, where $I$ is the defining ideal of the curve. The value semigroup $\Gamma$ associated to $\phi_k$ (or to the curve) is the set $$\Gamma= \{ord_y(\phi_k^*(h));\ h\in \O\setminus \{0\}\}.$$ If $\Gamma=\langle v_1,\ldots, v_r \rangle $, since $\phi_k$ is a primitive parameterization, then $gcd(v_1, \ldots, v_r)=1$ which implies that $\Gamma$ has a conductor $c$. If we denote by $\overline{\O}$ the integral closure of $\O$ in its field of fractions, then $\overline{\O}\simeq \Cp\{y\}$. The conductor ideal of $\O$ in $\overline{\O}$ is defined by $(\O:\overline{\O})=\{h \in \O;\ h\overline{\O} \subseteq \O\}$ which is an ideal of $\O$ satisfying that $(\O:\overline{\O})=\langle y^c \rangle$. Moreover, $(\O:\overline{\O})$ is also an ideal of $\overline{\O}$. In fact, the conductor ideal is the largest ideal of $\overline{\O}$ that is contained in $\O$.

 We consider the Ap\'{e}ry set $B=\{a_0,a_1, \ldots, a_{m_1-1}\}$ of $\Gamma$ with respect to the multiplicity $m_1$ of $\phi_k$. Both the semigroup and the set of gaps $L$ of $\Gamma$ can be determined in terms of the Ap\'{e}ry set $B$. In fact, any element $s\in \Gamma$ can be written as $s=a+\beta m_1$ for some $a\in B$ and $\beta \in \mathbb{N}\cup \{0\}$ and  \begin{equation}
\label{GapsAperySeq}
L=\left\{a_j-\beta_jm_1; \ j=1,\ldots, m_1-1, \ 1 \leq \beta_j \leq \left[\frac{a_j}{m_1}\right]\right\},
\end{equation}
where $[\cdot]$ denotes the integral part. The delta invariant of $\phi_k$ is given by $$\delta_{\phi_k}=dim_{\Cm}
\frac{\O_1}{\phi_k^*(\O_{k})}=dim_{\Cm}
\frac{\Cp\{y\}}{\phi_k^*(\O)}.$$  Thus, the previous dimension coincides with the number of elements in $\mathbb{N}$ that do not belong to the semigroup
$\Gamma$, that is, with the number of gaps of $\Gamma$. Thus, $\phi_k$ is a primitive parameterized curve if and only if $\delta_{\phi_k}$ is finite if and only if $\phi_k$ is $\A$-finite. For a parameterized plane curve $\phi_2$, the conductor $c$ of the semigroup is the Milnor number $\mu$ of the defining equation of the curve. Moreover, in this case
$$c=\mu = 2 \cdot \delta_{\phi_2}.$$
\begin{example}
	\label{Semigroup-4-5}
The semigroup associated to the parameterized curve $(y^4,y^5+y^6)$ is the set $\Gamma=\langle 4, 5 \rangle=\{0,4,5,8,9,10,12,13,14,\ldots\}$, the conductor of $\Gamma$ is $c=12$ and its set of gaps is $L=\{1,2,3,6,7,11\}$.
\end{example}

\end{itemize}

\section{Elementary join maps}
In this section we define an ``operation'' between two map-germs $F$ and $G$, called the join of $F$ and $G$, in order to produce a new map-germ that carries, in a certain sense, the $F$ and $G$ structures.

\begin{definition}
	Let $F:(\Cp^n,0)\to (\Cp^{p_1},0)$ and $G:(\Cp^m,0)\to (\Cp^{p_2},0)$ with $p_1 \geq n$ and $p_2 \geq m$. The \textit{join} of $F$ and $G$ (with respect to $H$), denoted by $(F*G)_{H}$, is a map from $(\Cp^{n+m},0)$ to $(\Cp^p,0)$ such that
	$$(F*G)_{H}(x,y)=(F(x),G(y),H(x,y)),$$ for some map-germ $H:(\Cp^{n+m},0)\to (\Cp^{p_3},0)$ in which $p=p_1+p_2+p_3$.
\end{definition}

Naturally we can extend the above definition for a finite number of maps $F_1,\ldots, F_r$ with respect to $H$ of the form: $(F_1*\ldots *F_r)_{H}(x_1,\ldots, x_r):=(F_1(x_1),\ldots, F_r(x_r),H(x_1,\ldots, x_r))$.

\vspace{0.2cm}
The $\A$-class of $(F*G)_{H}$ is well defined in the sense that if $F \underset{\A}{\sim} \widetilde{F}$ and $G \underset{\A}{\sim} \widetilde{G}$, then $(F*G)_{H} \underset{\A}{\sim} (\widetilde{F}*\widetilde{G})_{\widetilde{H}}$, for some map-germ $\widetilde{H}$ that is $\A$-equivalent to $H$.

\vspace{0.2cm}
Defining some special $\A$-finite map-germs that are given as a join of the identity and a parameterized curve or a join of $\A$-finite parameterized curves, we produce new families of $\A$-finite map-germs $f:(\Cp^n,0)\to (\Cp^p,0)$ with $p\geq 2n$.

In the next two results, we provide two families of $\A$-finite maps of corank $1$ and $n$, respectively, that allow us to characterize any $\A$-finite monomial map of corank $q$.

\begin{proposition}
	\label{FamilyA-Finite}
Let $\phi_k$ be the germ of a parameterized curve in $\Cp^k$ as in (\ref{ParamCurve}) and $I$ the identity in $\Cp^{n-1}$ with $n,k \geq 2$. Then the join $(I*\phi_k)_{H} :(\Cp^n,0)\to (\Cp^{2n+k-2},0)$, with respect to $H(\underline{x},y)=(x_1^{\lambda_1}y, \ldots, x_{n-1}^{\lambda_{n-1}}y)$, that is,
$$(I*\phi_k)_{H}(\underline{x},y):=(\underline{x},\ \phi_k(y),\ x_1^{\lambda_1}y, \ldots, x_{n-1}^{\lambda_{n-1}}y)$$
is an $\A$-finite map-germ of corank $1$, where $\underline{x}=(x_1, \ldots, x_{n-1})$ and $\lambda_i \geq 1$ for $i=1, \ldots, n-1$.
\end{proposition}
\Dem The proof that $f:=(I*\phi_k)_{H}$ is $\A$-finite follows a similar idea of \cite[Theorem 2.4]{RodRuas}. We can verify that $f$ is one-to-one, since $\phi_k$ is a primitive parameterization and in $\Cp^n\setminus \{0\}$ the map $df_x$ has maximal rank. Thus $f$ is stable outside the origin, and therefore, $\A$-finite by Mather-Gaffney geometric criterion. However we give an algebraic proof in order to estimate the delta invariant of $f$.

Let $\Gamma$ be the semigroup of $\phi_k$ with conductor $c$ and $L=\mathbb{N}\setminus \Gamma$ the set of gaps of $\Gamma$. We will prove that $\delta_f$ is finite since the following elements belong to $f^*(\O_p)$, where $p:=2n+k-2$:
\begin{enumerate}
	\item[$(a)$]  $x_1^{r_1}\cdots x_{n-1}^{r_{n-1}}y^s\in f^*(\O_p)$, for all $r_i\geq 0$ and $s=0$ or $s\geq c\ $;
	
	For each $i\in \{1,\ldots, n-1\}$:
	\item[$(b)$] $x_1^{r_1}\cdots x_i^{r_i}\cdots x_{n-1}^{r_{n-1}}y^s \in f^*(\O_p)$, with $s\in \Gamma$, $s<c\ $, $r_i\geq a\lambda_i$ for some $a$ depending on $s$ and for all $r_j\geq 0$ with $i\neq j\in \{1,\ldots, n-1\}$;
	
	\item[$(c)$] $x_1^{r_1}\cdots x_i^{r_i}\cdots x_{n-1}^{r_{n-1}}y^l \in f^*(\O_p)$, with $l\in L$, $r_i\geq (m_1-1)\lambda_i$ and for all $r_j\geq 0$ with $i\neq j\in \{1,\ldots, n-1\}$.
\end{enumerate}

In fact, it is clear that $x_1^{r_1}\cdots x_{n-1}^{r_{n-1}}\in f^*(\O_p)$ for all $r_i\geq 0$. On the other hand, by definition of the conductor of a semigroup, any integer $s\geq c$ satisfies $s\in \Gamma$ (see Section \ref{Sec-Notations}). We recall that the conductor ideal satisfies $(\O:\overline{\O})=\langle y^c \rangle$. Since $(\O:\overline{\O})$ is an ideal of $\O$ (and of $\overline{\O}$) we have for any $s\geq c$ that $y^s=y^{s-c}y^c \in (\O:\overline{\O})\subset \O$. Moreover $\O \simeq \Cp\{\phi_k(y)\}$, thus there exists $h\in \Cp\{X_1,\ldots, X_k\}$ such that $y^s=h(\phi_k) \in \phi_k^*(\O)\subset f^*(\O_p)$.

 For item $(b)$, if $s\in \Gamma$ with $s < c$, we only guarantee that there exists $h\in \Cp\{X_1,\ldots, X_k\}$ such that $h(\phi_k)=y^s+h.o.t. \in \phi_k^*(\O)$. But as we remarked before (\ref{GapsAperySeq}), we can write $s=a+\beta m_1$, for some $a$ in the Ap\'{e}ry set $B$ of $\Gamma$ and $\beta \in \mathbb{N}\cup \{0\}$. Thus for each $i=1, \ldots, n-1$ we have $$x_i^{a\lambda_i}y^s=x_i^{a\lambda_i}y^{a+\beta m_1}=(x_i^{\lambda_i}y)^{a}(y^{m_1})^{\beta} \in f^*(\O_p).$$ Therefore, $x_1^{r_1}\cdots x_i^{r_i}\cdots x_{n-1}^{r_{n-1}}y^s \in f^*(\O_p)$, for any $r_i\geq a\lambda_i$ and $r_j\geq 0$ with $i\neq j\in \{1,\ldots, n-1\}.$

In $(c)$, given $l\in L$ it follows by (\ref{GapsAperySeq}) that $l=a-\beta m_1$ for some $a\in B$ and $\beta \in \mathbb{N}$, where $1 \leq \beta \leq [\frac{a}{m_1}]$. Since $m_1[\frac{a}{m_1}] \leq a < m_1 ([\frac{a}{m_1}]+1)$ and denoting $d_{i1}=\lambda_{i}(a-m_1[\frac{a}{m_1}])$ for each $i=1,\ldots, m_1-1$ we have
\begin{equation}
\label{GapElemTh1}
x_i^{d_{i1}}y^{l}=(x_i^{\lambda_i}y)^{a-m_1[\frac{a}{m_1}]}(y^{m_1})^{[\frac{a}{m_1}]-\beta}\in f^*(\O_p).
\end{equation}
Thus $x_1^{r_1}\cdots x_i^{r_i}\cdots x_{n-1}^{r_{n-1}}y^l \in f^*(\O_p)$, for each $l\in L$ and $r_i\geq d_{i1}$, where $d_{i1}$ depends on $l \in L$. These computations are sufficient to guarantee that $\delta_f$ is finite. More strongly, one has that the power of $x_i$ does not depend on the gap. In fact, by the inequality $a +1 \leq m_1 ([\frac{a}{m_1}]+1)$ we have that $d_{i1}\leq (m_1-1)\lambda_i$ and therefore, by (\ref{GapElemTh1}), for any $l\in L$
$$x_i^{(m_1-1)\lambda_i}y^l=x_i^{(m_1-1)\lambda_i-d_{i1}}\cdot (x_i^{d_{i1}}y^l) \in f^*(\O_p).$$ \cqd

\begin{example}
	In this example we illustrate the previous result. Let $f:(\Cp^3,0)\to (\Cp^7,0)$ be an $\A$-finite map-germ given by $f(x_1,x_2,y)=(x_1,x_2,y^4,y^6+y^9,y^{11}+y^{18},x_1^2y,x_2y)$. The semigroup associated to the parameterized curve $(y^4,y^6+y^9,y^{11}+y^{18})$ is the set $\Gamma=\langle 4,6,11\rangle$, the conductor and the set of gaps of $\Gamma$ are $c=14$ and $L=\{1,2,3,5,7,9,13\}$, respectively. From the proof of Propositions \ref{FamilyA-Finite} and \ref{Estimating-delta} and some calculations we obtain that a basis for $\frac{\O_3}{f^*(\O_7)}$ is given by  $$\begin{array}{l}
	\gamma=	\{y^{l},x_1y^l,x_2y^2,x_2y^3,x_1^2y^2,x_1^2y^3, x_1^3y^2,x_1^3y^3,x_2^2y^3,x_1x_2y^2, \vspace{0.1cm}\\
	\hspace{1cm} x_1x_2y^3, x_1x_2^2y^3, x_1^2x_2y^3,x_1^3x_2y^3,x_1^4y^3,x_1^5y^3; \ \mbox{for all} \ \ l \in L\}.\end{array}$$ Thus $\delta_f=28$.	
\end{example}

In Proposition \ref{FamilyA-Finite}, we obtain a family of $\A$-finite map-germs joining the identity and one parameterized curve. A possible question is: Given two parameterized curves $\phi_1$ and $\phi_2$, is there a map $H$ such that the map-germ $(\phi_1 * \phi_2)_H$ is $\A$-finite? We illustrate this problem with an example to motivate the next result in which we join a finite number of parameterized curves.

\begin{example}
Let $g(x,y)=(x^3,x^4,y^5,y^6,x^2y,xy^3) \in  \O(2,6)$ be a corank $2$ map-germ. The semigroups of the parameterized curves $(x^3,x^4)$ and $(y^5,y^6)$ are respectively $\Gamma_1=\langle 3,4\rangle$ and $\Gamma_2=\langle 5,6\rangle$ with its corresponding set of gaps $L_1=\{1,2,5\}$ and $L_2=\{1,2,3, 4,7,8,9,13,14,19\}$. It is possible to verify that $\delta_g=48$, since a basis to $\frac{\O_2}{g^*(\O_6)}$ is the set $$\begin{array}{l}
\gamma=	\{x^{l_1},y^{l_2},xy,xy^2,xy^4,\ldots, xy^7,xy^{10},xy^{11},xy^{12}, xy^{16},xy^{17},xy^{22},x^2y^2,\ldots, x^2y^5,x^2y^8, \vspace{0.1cm}\\
	\hspace{1cm} x^2y^9,x^2y^{10},x^2y^{14},x^2y^{15},x^2y^{20},x^3y,x^3y^2,x^3y^3,x^3y^7,x^3y^8,x^3y^{13},x^4y,x^4y^4,x^5y^2, \vspace{0.1cm} \\
	\hspace{1cm} x^5y^4,x^6y^2, x^7y, x^9y^2; \ \mbox{for all} \ l_i \in L_i, \ i=1,2\}.\end{array}$$ 	
\end{example}
With the purpose of constructing an $\A$-finite family joining a finite number of curves, we will introduce some notation: for each $j=1,\ldots, n$ let $\phi_{k_j}:(\Cp,0)\to (\Cp^{k_j},0)$ be a primitive parameterized curve in $\Cp^{k_j}$, with multiplicity $m_{j} > 1$, given (up to $\A$-equivalence) by $$\phi_{k_j}(y_j)=(y_j^{m_j},y_j^{m_{j2}}+\sum_{i>m_{j2}}a_{i2}y_j^{i}, \ldots, y_j^{m_{jk_j}}+\sum_{i>m_{jk_j}}a_{ik_j}y_j^{i}),$$ with $1 < m_j:=m_{j1} \leq m_{j2}\leq \ldots \leq m_{jk_j}$. We denote by $H_j:(\Cp^n,0)\to (\Cp^{n-1},0)$ the map-germ
 \begin{equation}
 \label{ProdYiYj}
H_j(\underline{y}):=(y_j^{\mu_{j1}}y_1,\ldots,y_j^{\mu_{j(j-1)}}y_{j-1},y_j^{\mu_{j(j+1)}}y_{j+1},\ldots,y_j^{\mu_{jn}}y_{n}),
\end{equation} where $\underline{y}:=(y_1,\ldots, y_{n})$. For each $j\in \{1,\ldots, n\}$ the exponent $\mu_{js}$ of $y_j$ satisfy $\mu_{js}\geq 1$ for any $s\in \{1,\ldots, n\}$ with $s\neq j$.

\begin{proposition}
	\label{PropCorankn-map}
	With the previous notations, let $H:(\Cp^n,0)\to (\Cp^{n(n-1)},0)$ be a map-germ given by $H(\underline{y})=(H_1(\underline{y}),\ldots, H_{n}(\underline{y}))$ with $H_j$ as in (\ref{ProdYiYj}) and $n\geq 2$. Then $(\phi_{k_1}*\ldots*\phi_{k_n})_H:(\Cp^n,0) \to (\Cp^{p},0)$, that is,
	$$(\phi_{k_1}*\ldots*\phi_{k_n})_H(\underline{y})=(\phi_{k_1}(y_1), \ldots,\phi_{k_n}(y_n), H(\underline{y}))$$ is an $\A$-finite map-germ of corank $n$ in which $p=\sum_{j=1}^nk_j + n(n-1)$, $p \geq 2n$.
\end{proposition}
\Dem The Mather-Gaffney geometric criterion can be used, as in Proposition \ref{FamilyA-Finite}, to prove that $g:=(\phi_{k_1}*\ldots*\phi_{k_n})_H$ is $\A$-finite. However, we also present an algebraic proof in this case. Let $\Gamma_j$ be the semigroup of the parameterized curve $\phi_{k_j}$ with conductor $c_j$. We will verify that the following elements are in $g^*(\O_p)$:
\begin{enumerate}
	\item[$(a1)$] For each $s_i \geq c_i$, then $y_1^{s_1}\cdots y_i^{s_i}\cdots y_n^{s_n}\in g^*(\O_p) $ for any $s_j\geq c_j$ or $s_j= 0$ with $i, j\in \{1,\ldots, n\}$, $i\neq j$;
	
	\item[$(b1)$]  Given $s_i\in \Gamma_i$, $s_i < c_i$, $i\in \{1,\ldots, n\}$, then $y_t^{s_t}y_i^{s_i} \in g^*(\O_p)$ for any $t \in \{1,\ldots, n\}$, $t \neq i$, $s_t \geq c_t+a\cdot \mu_{ti}$ with $a$ depending on $s_i$;

	\item[$(c1)$] Given $l_i \in \mathbb{N}\setminus\Gamma_{i}$ with $i \in \{1,\ldots, n\}$, then $y_t^{s_t}y_{i}^{l_i} \in g^*(\O_p)$ for any $t\in \{1,\ldots, n\}$, $t\neq i$, $s_t \geq c_t+\alpha$ with $\alpha$ depending on $l_i$.
\end{enumerate}

The item $(a1)$ follows as in the proof of Proposition \ref{FamilyA-Finite} $(a)$. For $(b1)$ taking $s_i=a+\beta m_i\in \Gamma_i$ with $a$ in the Ap\'{e}ry sequence of $\Gamma_i$ and $\beta \in \mathbb{N}$ (see the proof of Proposition \ref{FamilyA-Finite} $(b)$), we have $$y_t^{a\cdot \mu_{ti}}y_i^{s_i}=(y_t^{\mu_{ti}}y_i)^a(y_i^{m_i})^{\beta}\in g^*(\O_p).$$
Since $y_t^{s}\in \phi_{k_t}^*(\O_{k_t})\subset g^*(\O_p)$ for all $s \geq c_t$, then for any $s_t \geq c_t+a\mu_{ti}$ it follows from the expression above that $y_t^{s_t}y_i^{s_i}=y_t^{s_t-a\mu_{ti}}\cdot (y_t^{a\mu_{ti}}y_i^{s_i}) \in g^*(\O_p)$.

In $(c1)$, given $l_i \in \mathbb{N}\setminus \Gamma_i$, as in (\ref{GapsAperySeq}) we can write $l_i=a-\beta m_i$, for some $a$ in the Ap\'{e}ry sequence of $\Gamma_i$ and $1 \leq \beta \leq [\frac{a}{m_{i}}]$. Then for any $t\in \{1,\ldots, n\}$ with $t\neq i$ and denoting $d_{ti}=\mu_{ti}(a-m_i[\frac{a}{m_i}])$ we obtain (as in (\ref{GapElemTh1})) that
$$y_t^{d_{ti}}y_i^{l_i}=(y_t^{\mu_{ti}}y_i)^{a-m_i[\frac{a}{m_i}]}(y_i^{m_i})^{[\frac{a}{m_i}]-\beta}\in g^*(\O_p).$$
Thus for any $s_t \geq c_t+d_{ti}$ the elements $y_t^{s_t}y_i^{l_i}=y_t^{s_t-d_{ti}}(y_t^{d_{ti}}y_i^{l_i})\in g^*(\O_p)$.

More generally, if we consider distinct $i_1,\ldots,i_q \in \{1,\ldots, n\}$, $s_{\rho}\in \Gamma_{i_{\rho}}$, with $s_{\rho}=0$ or $s_{\rho} < c_{i_{\rho}}$ (where $c_{i_{\rho}}$ is the conductor of $\Gamma_{i_{\rho}}$) for $\rho\in \{1, \ldots, q\}$, then we can write $s_{\rho}=a_{\rho}+\beta_{\rho}m_{i_{\rho}}$ with $a_{\rho}$, $\beta_{\rho}$ and $m_{i_{\rho}}$ as before. Thus, for any $t\in \{1,\ldots, n\}\setminus \{i_1,\ldots,i_q\}$ we have $$y_t^{(a_1\mu_{ti_1}+\ldots + a_q\mu_{ti_q})}y_{i_1}^{s_1}\cdots y_{i_q}^{s_q}\ = \ (y_t^{\mu_{ti_1}}y_{i_1})^{a_{1}}\cdots (y_t^{\mu_{ti_q}}y_{i_q})^{a_{q}}(y_{i_1}^{m_{i_1}})^{\beta_1}\cdots (y_{i_q}^{m_{i_q}})^{\beta_q}\in g^*(\O_p).$$
Then, for all $r_t \geq c_t+a_1\mu_{ti_1}+\ldots + a_q\mu_{ti_q}$ we obtain $y_t^{r_t}y_{i_1}^{s_1}\cdots y_{i_q}^{s_q} \in g^*(\O_p)$, and as consequence $y_1^{r_1}\cdots y_n^{r_n}(y_t^{r_t}y_{i_1}^{s_1}\cdots y_{i_q}^{s_q}) \in g^*(\O_p)$ for any $r_j \geq c_j$ or $r_j =0$ for all $j \in \{1,\ldots, n\}\setminus \{t,i_1,\ldots,i_q\}$.

\vspace{0.1cm}
Finally, given $l_{\rho}=a_{\rho}-\beta_{\rho} m_{i_{\rho}} \in \mathbb{N}\setminus\Gamma_{i_{\rho}}$,  $\rho=1,\ldots, q$, then for each $t\in \{1,\ldots, n\}\setminus \{i_1,\ldots,i_q\}$ and denoting $d_{ti_{\rho}}=\mu_{ti_{\rho}}(a_{\rho}-m_{i_{\rho}}[\frac{a_{\rho}}{m_{i_{\rho}}}])$, as in the proof of (c1), $y_t^{d_{ti_{\rho}}}y_{i_{\rho}}^{l_i}\in g^*(\O_p)$. Thus
$$y_t^{(d_{ti_1}+\ldots + d_{ti_q})}y_{i_1}^{l_1}\cdots y_{i_q}^{l_q}\ = \ (y_t^{d_{ti_1}}y_{i_1}^{l_{1}})\cdots (y_t^{d_{ti_q}}y_{i_q}^{l_{q}})\in g^*(\O_p).$$
Therefore, for all $r_t \geq c_t+d_{ti_1}+\ldots + d_{ti_q}$ and for all $j \in \{1,\ldots, n\}\setminus \{t,i_1,\ldots,i_q\}$ satisfying $r_j \geq c_j$ or $r_j=0$, the elements $y_1^{r_1}\cdots y_n^{r_n}\cdot (y_t^{r_t}y_{i_1}^{l_1}\cdots y_{i_q}^{l_q})\in g^*(\O_p)$.

A general element in $g^*(\O_p)$ can be obtained as a combination of the previous elements. Thus, $\delta_g$ is finite.\cqd

\begin{definition}
	The maps $(I*\phi_{k})_H$ and $(\phi_{k_1}*\ldots*\phi_{k_n})_H$ as in Proposition \ref{FamilyA-Finite} and Proposition \ref{PropCorankn-map}, respectively, are called \textit{elementary joins}.
\end{definition}

\begin{remark}
If $f:=(I*\phi_k)_H$ and $g:=(\phi_{k_1}*\ldots*\phi_{k_n})_H$ are the elementary join maps, and $p$ denotes its corresponding target dimension, then we can produce new $\A$-finite map-germs, for instance, $F(\underline{x},y)=(f(\underline{x},y),h(\underline{x},y))$ and $G(\underline{y})=(g(\underline{y}),k(\underline{y}))$, where $h,k \in \O_n$. 
\end{remark}

\begin{example}
	Let $f_a(x_1,x_2,y)=(x_1,x_2,y^4,y^5,x_1y,x_2y,ax_1x_2y^3) \in  \O(3,7)$ be a corank $1$ map-germ, where $a=0$ or if $a\neq 0$ then, up to $\A$-equivalence, we take $a=1$. In both cases, $f_a$ is $\A$-finite. The semigroup of the parameterized curve $(y^4,y^5)$ is $\Gamma=\langle 4, 5 \rangle$ and its set of gaps is $L=\{1,2,3,6,7,11\}$ (see Example \ref{Semigroup-4-5}). Notice that a basis for $\frac{\O_3}{f_a^*(\O_7)}$ is $$\beta= \{y^l,x_iy^2,x_iy^3,x_iy^7, x_1x_2y^3, x_i^2y^3; \ i=1,2, \forall \ l \in L\}$$ if $a=0$ and $\beta \setminus \{x_1x_2y^3\}$ if $a=1$. Thus, $\delta_{f_0}=15$ and $\delta_{f_1}=14$.
\end{example}

\section{General form of $\A$-finite monomial maps of corank $q$}

Our purpose in this section is to give a general form for any map-germ $f:(\Cp^n,0) \to (\Cp^p,0)$ with $p\geq 2n$, satisfying the properties that $f$ is monomial and $\A$-finitely determined. Moreover we obtain normal forms of some classes of map-germs.

In the next lemma we prove that the class of parameterized curves appears naturally in the characterization of $\A$-finite monomial maps in $\O(n,p)$ for $p\geq 2n$.

\begin{lemma}
	\label{RemarkProjCurve}
	If $g:(\Cp^n,0) \to (\Cp^{p},0)$ with $p \geq 2n$ is a singular $\A$-finite monomial map-germ, then there exists at least one $\A$-finite parameterized monomial curve $\phi_k:(\Cp,0)\to (\Cp^{k},0)$ with $2 \leq k<p$ such that $g$ is $\A$-equivalent to $f(\underline{x},y)=\left(f_1(\underline{x},y),\ldots,f_{p-(n+k-1)}(\underline{x},y), x_1^{\lambda_1}y,\ldots, x_{n-1}^{\lambda_{n-1}}y ,\phi_k(y)\right)$, where $\underline{x}=(x_1,\ldots, x_{n-1})$.
\end{lemma}
\Dem Denote by $(x_1,\ldots, x_{n})$ the variables of $g=(g_1,\ldots, g_p)$. By the hypothesis that $\delta_g$ is finite, we have that $\M_n^{\beta} \subset g^*(\O_p)$ for some $\beta \in \mathbb{N}$. In particular, $x_1^{\gamma_1}, \ldots, x_n^{\gamma_n} \in g^*(\O_p)$ for all $\gamma_i \geq \beta$, $i=1,\ldots, n$. Since $g$ is monomial, it follows that for each $i \in \{1,\ldots, n\}$ there exists a component of $g$ given by $x_i^{r_i}$ for some $r_i \geq 1$. Moreover, being $g$ singular, then at least one variable, namely $x_n$, is such that $r_n\geq 2$, or else $g$ is $\A$-equivalent to an immersion. Let us rename $x_n$ and $r_n$ respectively by $y$ and $m_1$.
	
In this way, in order to obtain all powers of the variable $y$ of the form $y^b \in g^*(\O_p)$ for all $b \geq \beta$ with $g$ monomial, it is necessary to have other components of $g$ given as powers of $y$ so that there is an $\A$-finite parameterized curve in $y$. Then there exists $k$ with $1 < k< p$ such that $g$ can be given by $g(\underline{x},y)=(g_1(\underline{x},y),\ldots, g_{p-k}(\underline{x},y),y^{m_1}, \ldots, y^{m_k})$ with $\underline{x}=(x_1,\ldots, x_{n-1})$ and $g_i(\underline{x},0)=0$ for $i=1,\ldots, p-k$ and $gcd(m_1, \ldots, m_k)=1$. If we suppose $1< m_1 < m_2 < \ldots < m_k$, by the same argument that $g$ is monomial and $\delta_g$ is finite, the elements $x_1^{\lambda_1}y,  \ldots, x_{n-1}^{\lambda_{n-1}}y$ appear as components of $g$, for some $\lambda_i \in \mathbb{N}$, $i=1,\ldots, n-1$. Otherwise we could not obtain these elements in terms of other monomial components of $g$, since the multiplicity $m_1 >1$. Then, $g$ is $\A$-equivalent to $f(\underline{x},y)=(f_1(\underline{x},y),\ldots,f_{p-(n+k-1)}(\underline{x},y), x_1^{\lambda_1}y,\ldots, x_{n-1}^{\lambda_{n-1}}y,\phi_k(y))$.
\cqd

In Theorem \ref{ClassMonomProjCurve} and Corollaries \ref{ClassifMon-n2n} and \ref{ClassMon-General} we characterize some classes of $\A$-finite monomial map-germs in terms of elementary joins.

\begin{theorem}
	\label{ClassMonomProjCurve} Let $\phi_k(y)=(y^{m_1}, \ldots, y^{m_k})$ be a parameterized monomial curve in $\Cp^k$, with $1< m_1 < m_2< \ldots < m_k$, $gcd(m_1,\ldots, m_k)=1$ and $n,k \geq 2$. Then, a singular monomial map-germ $f:(\Cp^n,0) \to (\Cp^{2n+k-2},0)$ satisfying $f(\underline{0},y)=(0,\ldots, 0,\phi_k(y))$ is $\A$-finite if and only if $f$ is $\A$-equivalent to
	 $$(I*\phi_k)_H(\underline{x},y)=(\underline{x},y^{m_1}, \ldots, y^{m_k}, x_1^{\lambda_1}y, \ldots, x_{n-1}^{\lambda_{n-1}}y),$$ where $\underline{x}=(x_1, \ldots, x_{n-1})$, $I$ is the identity in $\Cp^{n-1}$ and $\lambda_i \geq 1$ for $i=1, \ldots, n-1$.
\end{theorem}
\Dem Notice that the sufficient condition follows from Proposition \ref{FamilyA-Finite}.

By hypothesis $f(\underline{x},y)=(f_1(\underline{x},y),\ldots, f_{2n-2}(\underline{x},y),\phi_k(y))$ with $f_i(\underline{0},y)=0$ for $i=1,\ldots, 2n-2$. From the proof of Lemma \ref{RemarkProjCurve}, there exist $n-1$ components of $f$ of the form $x_i^{r_i}$ with $r_i\geq 1$ for $i=1, \ldots, n-1$. We claim that $r_i=1$, for all $i=1,\ldots, n-1$. In fact, if we assume by contradiction that some $r_i$ is greater than $1$, say $r_1$, we have an $\A$-finite parameterized curve in the variable $x_1$ as component of $f$. Suppose, without loss of generality, that such curve is given by $(x_1^{n_1},x_1^{n_2})$  with $r_1:=n_1$, $1 < n_1 < n_2$ and $gcd(n_1,n_2)=1$. Then, by Lemma \ref{RemarkProjCurve}, the elements: $x_1^{n_1},x_1^{n_2}, x_2^{s_2}x_1, \ldots, x_{n-1}^{s_{n-1}}x_1, y^{s_1}x_1,y^{m_1},\ldots, y^{m_k},x_1^{\lambda_1}y,  \ldots, x_{n-1}^{\lambda_{n-1}}y$ are components of $f$, for some $\lambda_i,s_i\in {\mathbb N}$, $i=1,\dots, n-1$. But this is not possible, since the target dimension is $p=2n+k-2$. Then $r_1=\ldots=r_{n-1}=1$ and we conclude that $f$ is a corank one map. Therefore, $f$ is $\A$-equivalent to a join of the identity $I$ in $\Cp^{n-1}$ and $\phi_k$, with respect to the map $H:(\Cp^n,0) \to (\Cp^{n-1},0)$, $H(\underline{x},y)=(x_1^{\lambda_1}y, \ldots, x_{n-1}^{\lambda_{n-1}}y)$. \cqd

It is interesting to observe that if $f \in \M_n\O(n,2n+k-2)$ is a singular monomial and $\A$-finite map-germ, then $f$ cannot have as component a parameterized curve $\phi_j:(\Cp,0)\to (\Cp^j,0)$ for $j> k$, by a target dimension argument.

As a consequence of Lemma \ref{RemarkProjCurve} and Theorem \ref{ClassMonomProjCurve}, we extend for maps in $\M_n\O(n,2n)$ the classification of parameterized monomial surfaces in $\Cp^4$ obtained in \cite[Theorem 2.4]{RodRuas}.

\begin{corollary}
	\label{ClassifMon-n2n}
Let $f:(\Cp^n,0) \to (\Cp^{2n},0)$ be a monomial map-germ with $n\geq 2$. Then, $f$ is $\A$-finite if and only if $f$ is $\A$-equivalent to
 \begin{enumerate}
 	\item[(i)] $(\underline{x},y,0,\ldots, 0)$ (immersion);
 	\item[(ii)] $(\underline{x},y^{m_1}, y^{m_2}, x_1^{\lambda_1}y, \ldots, x_{n-1}^{\lambda_{n-1}}y),$ where $\lambda_i \geq 1$ for $i=1, \ldots, n-1$ and $gcd(m_1,m_2)=1$.
 \end{enumerate}
 \end{corollary}

\vspace{0.2cm}
In what follows, we give a general form of an $\A$-finite monomial map in $\O(n,p)$ of corank $q$, for $p \geq 2n$. For each $j=1,\ldots, q$ the map $\phi_{k_j}(y_j)=(y_j^{m_j},y_j^{m_{j2}}, \ldots, y_j^{m_{jk_j}})$ denotes a primitive parameterized curve in $\Cp^{k_j}$ with $1 < m_j:=m_{j1} < m_{j2}< \ldots < m_{jk_j}$ and if $\underline{x}:=(x_1,\ldots, x_{n-q})$ and $\underline{y}:=(y_1,\ldots, y_{q})$, then $\underline{x}^{\lambda_{ij}}y_j:=(x_1^{\lambda_{1j}}y_j,\ldots, x_{n-q}^{\lambda_{(n-q)j}}y_j)$ and (as in (\ref{ProdYiYj})) $$H_j(\underline{y}):=(y_j^{\mu_{j1}}y_1,\ldots,y_j^{\mu_{j(j-1)}}y_{j-1},y_j^{\mu_{j(j+1)}}y_{j+1},\ldots,y_j^{\mu_{jq}}y_{q}),$$ where $\lambda_{ij} \geq 1$, $\mu_{jt}\geq 1$ for $i\in \{1, \ldots, n-q\}$, $j,t\in \{1, \ldots, q\}$ with $t\neq j$ and $I$ denotes the identity map in $\Cp^{n-q}$.

\begin{theorem}{(General form of $\A$-finite monomial maps)}
	\label{ClassifMonoCurve-Corankq}
	Let $f:(\Cp^n,0) \to (\Cp^{p},0)$ be a singular corank $q$ monomial map-germ with $p\geq 2n$. Then $f$ is $\A$-finite if and only if $f$ is $\A$-equivalent to
	\begin{equation}
	\label{GeneralFormCorankq}
	(I*\phi_{k_1}* \ldots * \phi_{k_q})_H(\underline{x},\underline{y}) = (I(\underline{x}),\phi_{k_1}(y_1), \ldots, \phi_{k_q}(y_q),H(\underline{x},\underline{y})),
	\end{equation}
	where $H(\underline{x},\underline{y})=(\underline{x}^{\lambda_{i1}}y_1, \ldots,\underline{x}^{\lambda_{iq}}y_q,H_1(\underline{y}),\ldots,H_q(\underline{y}),h(\underline{x},\underline{y}))$ for some monomial map $h(\underline{x},\underline{y})$ such that $h$ is zero on each coordinate axis. In particular, if all $\mu_{jt}= 1$ the map $H$ can given by $$H(\underline{x},\underline{y})=(\underline{x}^{\lambda_{i1}}y_1, \ldots, \underline{x}^{\lambda_{iq}}y_q, y_1y_2,\ldots, y_1y_q,y_2y_3,\ldots, y_2y_q,\ldots, y_{q-1}y_q,h(\underline{x},\underline{y})).$$
\end{theorem}
\Dem It follows from the proofs of Propositions \ref{FamilyA-Finite} and \ref{PropCorankn-map} that the map (\ref{GeneralFormCorankq}) is $\A$-finite.

Conversely, if $f$ is a corank $q$ map-germ, up to change of coordinates, we can assume that $f(\underline{x},\underline{y})=(\underline{x},f_{n-q+1}(\underline{x},\underline{y}),\ldots,f_p(\underline{x},\underline{y}))$. Since $f$ is $\A$-finite and monomial, it follows from the proof of Lemma \ref{RemarkProjCurve}, that for each $j=1,\ldots, q$ there exists $k_j$ with $1 < k_j \leq p-(n-q)$ such that $f$ has as components an $\A$-finite parameterized curve in $\Cp^{k_j}$ in the variable $y_j$, that is, $\phi_{k_j}(y_j)=(y_j^{m_j},y_j^{m_{j2}}, \ldots, y_j^{m_{jk_j}})$ where we suppose $1 <m_j <m_{j2}< \ldots < m_{jk_j}$. If $m_j=1$ for some $j=1, \ldots, q$, then the corank of $f$ would be less than $q$. Thus, $f$ is $\A$-equivalent to $(I*\phi_{k_1}* \ldots * \phi_{k_q})_H$ for some map-germ $H$. Similarly as in the proof of Theorem \ref{ClassMonomProjCurve}, the hypotheses that $f$ is monomial and $\delta_f$ finite imply that $H$ has to be given by $H(\underline{x},\underline{y})=(\underline{x}^{\lambda_{i1}}y_1, \ldots, \underline{x}^{\lambda_{iq}}y_q,H_1(\underline{y}),\ldots,H_q(\underline{y}),h(\underline{x},\underline{y}))$, for some monomial map-germ $h$ that, up to $\A$-equivalence, is zero on each coordinate axis. \cqd

\begin{remark}
Notice that in Theorem \ref{ClassifMonoCurve-Corankq}, the components of $h$ can be zero. Furthermore, if $\phi_{k_j}$ or $h$ are not monomials, we still have that the map (\ref{GeneralFormCorankq}) is $\A$-finite, that is, combining elementary joins we produce a lot of $\A$-finite maps.
\end{remark}

Interesting particular cases in Theorem \ref{ClassifMonoCurve-Corankq} are when $f$ is one of the elementary join maps:
$$\begin{array}{l}
\bullet \ corank\ 1 : \ (\underline{x}, y^{m_1}, \ldots,y^{m_r},x_1^{\lambda_1}y, \ldots, x_{n-1}^{\lambda_{n-1}}y);\vspace{0.2cm}\\
\bullet \ corank \ n : \ (\phi_{k_1}(y_1), \ldots,\phi_{k_n}(y_n),\ H_1(\underline{y}),\ldots, H_n(\underline{y})).\end{array}$$

\begin{corollary}
	\label{ClassMon-General}
	Let $f:(\Cp^n,0) \to (\Cp^{p},0)$ be a singular corank $q$ monomial map-germ, with $p\geq 2n.$ If $f$ is $\A$-finite, then $p\geq n(q+1)-q(q-1)/2$. In particular, if $p=n(q+1)-q(q-1)/2$, then $f$ is $\A$-finite if and only if $$f \ \underset{\A}{\sim} \ (I*\phi_{k_1}* \ldots * \phi_{k_q})_H,$$
	where for each $j=1,\ldots, q$ the map $\phi_{k_j}(y_j)=(y_j^{m_j},y_j^{m_{j1}})$ is a primitive parameterized plane curve and $H(\underline{x},\underline{y})=(\underline{x}^{\lambda_{i1}}y_1, \ldots, \underline{x}^{\lambda_{iq}}y_q, y_1y_2,\ldots, y_1y_q,y_2y_3,\ldots, y_2y_q,\ldots, y_{q-1}y_q)$.
\end{corollary}

Notice that by the previous result, in order to have an $\A$-finite monomial map-germ in $\M_n\O(n,p)$ of corank $n$, then $p \geq \frac{n(n+3)}{2}$. But for corank $1$ maps, we need $p\geq 2n$.

\begin{example}
 Let $h(x,y,z)=(x,y^2,y^3,z^2,z^3,xy,xz,yz)\in \O(3,8)$ be a corank $2$ map-germ. In this case, a basis for $\frac{\O_3}{h^*(\O_8)}$ is $\{y,z,y^2z,yz^2\}$ and thus $\delta_h=4$.
\end{example}

\begin{remark}
In Theorem \ref{ClassifMonoCurve-Corankq}, we obtain a general form of an $\A$-finite monomial map-germ from $(\Cp^n,0)$ to $(\Cp^{p},0)$ for $p \geq 2n$. In this sense, an interesting question is whether we can describe all $\A$-finite monomial map-germs $F:(\Cp^n,0)\to (\Cp^{n+s},0)$ with $s=1, \ldots, n-1$? The only classes we got in $\O(n,n+s)$ satisfying that the germ is $\A$-finite and monomial were the immersion and
$$F(x_1,\ldots, x_{n-1},y)=(x_1,\ldots, x_{n-1},x_{i_1}y,x_{i_2}y,\ldots, x_{i_s}y, y^{2}),\ $$ for each $s=1,\ldots, n-1$ and distinct $i_1,\ldots i_s \in \{1, \ldots, n-1\}$. Notice that $F$ is a stable map, that is, $\A_e cod(F)=0$. Moreover, $F$ is a generalization of the cross-cap parameterized by $f(x,y)=(x,xy,y^2)$ (see \cite{Whitney}).
\end{remark}

\section{Estimating some analytical invariants}

In this section, given $f$ in $\O(n,p)$ with $p\geq 2n$ we estimate the $\A_e$-codimension and the delta invariant of $f$, for the classes of map-germs of the previous sections. When $p=2n$, we obtain some formulas for the number of double points in a stable perturbation of $f$. In the last subsection, we study these invariants under some projection of $f$. This approach can be useful when the projection has good properties, such as being stable.

If $f\in \M_n\O(n,p)$ with $p\geq 2n$ satisfies $\delta_f < \infty$, then $f$ is $\A$-finite. Moreover, $$\Acod \leq {\mathcal L}_ecod(f)= p \cdot \delta_f.$$ In this way, it is interesting to obtain some estimates for $\delta_f$.

Let $F:(\Cp^n,0)\to (\Cp^{p_1},0)$, $G:(\Cp^m,0)\to (\Cp^{p_2},0)$ and $H:(\Cp^{n+m},0)\to (\Cp^{p_3},0)$ with $p_1 \geq n$, $p_2 \geq m$ such that $K=(F*G)_{H}:(\Cp^{n+m},0)\to (\Cp^p,0)$ is a join of $F$ and $G$ with respect to $H$, that is,
$K(x,y)=(F(x),G(y),H(x,y))$ where $p=p_1+p_2+p_3$.

Suppose $p_1\geq 2n$, $p_2 \geq 2m$, $\delta_F$ and $\delta_G$ finite. If $H$ is zero on each coordinate axis, then it is easy to see that
\begin{equation}
\label{DeltaJoinMaps}
\delta_K \geq \delta_F + \delta_G.
\end{equation}

\vspace{0.3cm}
For the next result, let $\phi_k(y)=(y^{m_1}, y^{m_2}+\sum_{i>m_2}a_{i2}y^{i},\ldots, y^{m_k}+\sum_{i>m_k}a_{ik}y^{i})$ be an $\A$-finite parameterized curve and consider the following notations: the integers $\lambda_{\min}:=\underset{1\leq i \leq n-1}{\min}\{\lambda_i\}$,  $\ \lambda:=(m_1-1)\sum_{i=1}^{n-1}\lambda_i-n+1$,
$$\begin{array}{l}
\kappa_1:=\displaystyle\sum_{i=1}^{n-1}\lambda_i +\scriptsize{\sum_{j=1}^{\lambda_{\min}}\left(\hspace{-0.15cm}\begin{array}{c} n+j-2 \\
	j \\ \end{array}\hspace{-0.15cm}\right)} - (n-1)\lambda_{\min}-n+2,\\ \kappa_2:=\displaystyle\sum_{j=1}^{\lambda}\scriptsize{\left(\hspace{-0.15cm}\begin{array}{c}
	n+j-2 \\
	j \\ \end{array}\hspace{-0.15cm}\right)}+1 \ \ \ and \ \ \ \kappa_3:=\displaystyle\sum_{i=1}^{n-1}\lambda_i -n.\end{array}$$

\begin{proposition}
	\label{Estimating-delta}
	Let $\phi_k$ be the germ of a parameterized curve in $\Cp^k$ as above with $n,k \geq 2$. Consider $f :(\Cp^n,0)\to (\Cp^{2n+k-2},0)$ the map-germ given by $f(\underline{x},y)=(\underline{x}, \phi_k(y), x_1^{\lambda_1}y, \ldots, x_{n-1}^{\lambda_{n-1}}y)$, where $\underline{x}=(x_1, \ldots, x_{n-1})$ and $\lambda_i \geq 1$ for $i=1, \ldots, n-1$. Then
	\begin{enumerate}
		\item[(i)] $\kappa_1 \delta_{\phi_k}  \leq \delta_{f}\leq  \kappa_2\delta_{\phi_k},$
		\item[(ii)] $(n-1)\kappa_3 \delta_{\phi_k}+k(m_1-2)(\kappa_3+2)\leq \A_ecod(f) \leq (2n+k-2)\kappa_2\delta_{\phi_k}$.
	\end{enumerate}
\end{proposition}
\Dem Denote $p=2n+k-2$ and $c$ the conductor of the semigroup $\Gamma$ of $\phi_k$. It follows from the proof of Proposition \ref{FamilyA-Finite} that $x_i^{r_i},y^b,x_i^{p_i}y^s, x_i^{q_i}y^l \in f^*(\O_p)$, for all $r_i \in \mathbb{N}$, $b\geq c$, $s\in \Gamma$ for some $p_i$ and $l \in L=\mathbb{N}\setminus \Gamma$ and $q_i \geq (m_1-1)\lambda_i$, with $i=1, \ldots, n-1$. Now we will analyse the element $y^s$ with $s\in \Gamma$ and $s$ less than the conductor $c$. If $s\in \Gamma$, then there exists $h \in \O$ such that $h\circ \phi_k(y)=y^s+a_{s+1}y^{s+1}+\ldots \in \phi_k^*(\O_k)\subset f^*(\O_p)$ with $a_j \in \Cp$ for $j=s+1, \ldots$. For each $j\in \Gamma\setminus \{0\}$ with $j>s$ in the expansion $h\circ \phi_k(y)$, there exists $h_j \in \O$ such that $h_j\circ \phi_k(y)=b_jy^j+\ldots$, with $b_j\neq 0$. Thus $$\overline{h}=h-\frac{a_j}{b_j}\sum_{j\in \Gamma} h_j \in \O$$ satisfy
$\overline{h}\circ \phi_k(y)=y^s+a_{l_1}y^{l_1}+\ldots +a_{l_t}y^{l_t}$ for some finite $l_1, \ldots, l_t \in L$ and $a_{l_i}\in \Cp$, since $y^b \in f^*(\O_p)$ for $b\geq c$. Denoting $\beta$ a basis of $\O_n/f^*(\O_p)$, we have that $y^l \in \beta$ for all $l\in L$. In this way, by the expression of $\overline{h}\circ \phi_k$, we do not consider in $\beta$ the element $y^s$, with $s\in \Gamma$. With a similar analysis $x_1^{r_1}\ldots x_{n-1}^{r_{n-1}}y^s \notin \beta$ for $r_i \geq 0$ and $s\in \Gamma$.

Thus $\beta$ is contained in $\{y^l, x_1^{r_1}\ldots x_{n-1}^{r_{n-1}}y^l; \ \forall \ l\in L, 1\leq d \leq \lambda\}$ with $d=\sum_{i=1}^{n-1}r_i$ and $\lambda:=(m_1-1)\sum_{i=1}^{n-1}\lambda_i-n+1.$ Therefore $$\delta_f \leq \delta_{\phi_k}+\sharp\{x_1^{r_1}\ldots x_{n-1}^{r_{n-1}}y^l; \ \forall \ l\in L, 1\leq d \leq \lambda\}=\kappa_2\delta_{\phi_k}.$$

 On the other hand, the elements $y^l, x_i^{r_i}y^l \in \beta$ for all $1 \leq r_i \leq \lambda_i -1$, $l \in L$ and $i=1, \ldots, n-1$. Moreover $x_1^{r_1}\ldots x_{n-1}^{r_{n-1}}y^l \in \beta$ for all $l \in L$ such that $1 \leq d \leq \lambda_{\min}$, eliminating the monomials $x_1^dy^l, \ldots, x_{n-1}^dy^l$ that have been considered in $\beta$. Thus, $$\delta_f \geq \sharp\{y^l, x_i^{r_i}y^l;\ \forall \ l\in L, 1 \leq r_i \leq \lambda_i -1\}+\sharp\{x_1^{r_1}\ldots x_{n-1}^{r_{n-1}}y^l;\forall \ l\in L, 1\leq d \leq \lambda_{\min}\}-(n-1)\lambda_{\min}\delta_{\phi_k}$$ and we obtain $\delta_f \geq \kappa_1\delta_{\phi_k}$.

The upper bound of $\A_ecod(f)$ is obtained from inequality $\Acod \leq p \delta_f$ and item $(i)$. The extended tangent space to the $\A$-orbit of $f$ is given by
{\scriptsize $$T\A_e(f)=\left\{\left(\begin{array}{cccc}
1 & \ldots & 0 & 0   \\
\vdots  & \ddots & \vdots & \vdots   \\
0 & \ldots & 1 & 0   \\
    0   & \ldots &   0    & m_1 y^{m_1-1} \\
       \vdots  & \cdots & \vdots & \vdots   \\
     0   & \ldots &   0    & m_k y^{m_k-1}+\ldots \\
 \lambda_1 x_1^{\lambda_1-1}y & \ldots & 0 & x_1^{\lambda_1}\\
  \vdots  & \cdots & \vdots & \vdots   \\
 0  & \ldots &  \lambda_{n-1} x_{n-1}^{\lambda_{n-1}-1}y & x_{n-1}^{\lambda_{n-1}}\\
 \end{array}
\right)\left(\begin{array}{c}
\epsilon_1(\underline{x},y)\\
\vdots \\
\epsilon_n(\underline{x},y)
\end{array}\right)+\left(\begin{array}{c}
\eta_1(f)\\
\vdots \\
\eta_{p}(f)\\
\end{array}\right);\begin{array}{c}
\epsilon_i \in \O_n, \\
\eta_j \in \O_{p},\\
i=1, \ldots, n, \\
j=1, \ldots, p \\
\end{array}\right\}$$}
We denote by $e_i$ the $i$-th element of the canonical basis of $\theta(p)$. For the lower bound of the $\A_e$-codimension of $f$, note that the following sets are not contained in $T\A_e(f)$:
$$\{\{y^l,x_i^{r_i}y^l\}e_{j};\ i\in \{1,\ldots, n-1\}, j \in \{n,\ldots, n+k-1\},  1\leq r_i \leq \lambda_i-1 \ \mbox{and}\ 1 \leq l \leq m_1-2 \}$$ and
$$\{\{x_i^{r_i}y^l\}e_{n+k+j-1};\ i,j\in \{1,\ldots, n-1\}, \forall \ l\in L, 1\leq r_i \leq \lambda_i-2 \ \mbox{if}\ i=j \ \mbox{and} \ 1 \leq r_i \leq \lambda_i-1\ \mbox{if}\ i\neq j \}.$$  
\cqd

For the map $(\phi_{k_1}*\ldots*\phi_{k_n})_H$ in Proposition \ref{PropCorankn-map} the estimate for the delta invariant is more general than the previous result. If $\Gamma_j$ is the semigroup of $\phi_{k_j}$ with conductor $c_j$, we denote
$c=\sum_{j=1}^nc_j$.

\begin{proposition}
	Let $g:(\Cp^n,0) \to (\Cp^{p},0)$ be the map-germ  as in Proposition \ref{PropCorankn-map} given by $$g(\underline{y})=(\phi_{k_1}(y_1), \ldots,\phi_{k_n}(y_n),\ H_1(\underline{y}),\ldots,H_n(\underline{y})),$$ where $H_j(\underline{y}):=(y_j^{\mu_{j1}}y_1,\ldots,y_j^{\mu_{j(j-1)}}y_{j-1},y_j^{\mu_{j(j+1)}}y_{j+1},\ldots,y_j^{\mu_{jn}}y_{n}),$
with $\underline{y}:=(y_1,\ldots, y_{n})$ and for each $j\in \{1,\ldots, n\}$ we have $t\in \{1,\ldots, n\}$ with $t\neq j$ and $\mu_{jt}\geq 1$. Then
	$$\sum_{j=1}^n\delta_{\phi_{k_j}}\ \leq \ \delta_g \ \leq \ \sum_{j=1}^c \left(\hspace{-0.15cm}\begin{array}{c} n+j-1 \\
	j \\ \end{array}\hspace{-0.15cm}\right).$$
\end{proposition}
\Dem It follows from (\ref{DeltaJoinMaps}) the lower estimate for $\delta_g$. On the other hand, if $\gamma$ denotes a basis for $\O_n/g^*(\O_p)$, by the proof of Proposition \ref{PropCorankn-map}, $\gamma \subset \{y_1^{s_1}\cdots y_n^{s_n};\ 1 \leq d \leq c, d=\sum_{j=1}^ns_j\}$. Thus, we obtain the result.\cqd

\vspace{0.1cm}
A map-germ $f$ is called $\A$-stable if $\theta(f)=\Tg$. Similarly as in \cite[Proposition 2.1]{RodRuas}, it is possible to prove for $p \geq 2n$ that $f$ is $\A$-stable if and only if $\delta_f=0$.

Let $f:(\Cp^n,0)\to (\Cp^{2n},0)$ be an $\A$-finite map. We recall that a stable perturbation $f_t:U \to (\Cp^{2n},0)$ of $f$, where $U$ is a small neighborhood of the origin, has a finite number $d(f)$ of transverse double points (see \cite{GolubitskyGuil}).

If $f$ is an $\A$-finite corank
one map-germ, it is well known that up to change of coordinates $f$ is given by
$f(\underline{x},y)= (\underline{x},f_n(\underline{x},y), \ldots,
f_{2n}(\underline{x},y))$, where $\underline{x}:=(x_1, \ldots, x_{n-1})$. In this case, the number of
double points of $f$ is
\begin{equation}
\label{d(f)G} d(f):=\displaystyle\frac{1}{2} dim_{\Cp}\
\frac{\O_{n+1}}{\varphi^*(\M_{n+1})\O_{n+1}}, \end{equation} where
$ \varphi: (\Cp^{n+1},0)  \to  (\Cp^{n+1},0)$, $(\underline{x},y,z) \mapsto
(\varphi_n,\ldots,\varphi_{2n})$ with $\varphi_i(\underline{x},y,z):=
\frac{f_i(\underline{x},z)-f_i(\underline{x},y)}{z - y}$, for each $i=n
\ldots, 2n$ (see \cite{Klotz}). Moreover, if $f$ is quasihomogeneous of type $(\omega_1, \ldots, \omega_n,d_n,\ldots, d_{2n})$, where $\omega_1, \ldots, \omega_{n-1}, \omega_n$ are the weights of the variables $x_1, \ldots, x_{n-1},y$ and $d_n, \ldots,d_{2n}$ are the weighted degrees of the components functions $f_n, \ldots, f_{2n}$ of $f$, respectively, then
\begin{equation}
\label{FormDoubPointsQuasih}
d(f)=\frac{\prod_{i=n}^{2n}(d_i-\omega_n)}{2\omega_1\ldots \omega_{n-1}\omega_n^2}.
\end{equation}
The previous formula was obtained by applying Bezout's generalized formula to $\varphi$ with $\omega_n$ the weight of the variable $z$ (see \cite{Klotz}).
\begin{proposition}
	\label{PtDupCurv} Let $\phi_2(y)=(y^{m_1},y^{m_2})$ with $1 < m_1
	< m_2$, $gcd(m_1,m_2)=1$ and for each $i=1,2$ the map $h_i(\underline{x},y)=
	y^{m_i} + \sum_{j>m_i}a_{ij}(\underline{x})y^{j}$ satisfies  $a_{ij}(\underline{x})
	\in \M_{n-1}$.
	\begin{enumerate}
		\item If $f(\underline{x},y)=(\underline{x},
		\underline{x}y, h_1(\underline{x},y), h_2(\underline{x},y))$, then
		$d(f)=\frac{(m_1-1)(m_2-1)}{2}=\delta_{\phi_2}$.
		\item If $f(\underline{x},y)=(\underline{x}, y^{m_1}, y^{m_2},x_1^{\lambda_1}y, \ldots, x_{n-1}^{\lambda_{n-1}}y)$, then $d(f)=\left(\prod_{i=1}^{n-1}\lambda_i\right)\delta_{\phi_2},$ for  $\lambda_i \geq 1$.
	\end{enumerate}
\end{proposition}
\Dem By definition of $d(f)$ in (\ref{d(f)G}), $\varphi_i(\underline{x},y,z)=x_{i-n+1}$ for $i=n, \ldots, 2n-2$, and denoting $P_{k}(z,y)=\frac{z^{k+1}-y^{k+1}}{z-y}$ we have
$$\varphi_{2n+i-2}(\underline{x},y,z)= P_{m_i-1}(z,y)+
\sum_{j> m_i}a_{ij}(\underline{x})P_{j-1}(z,y)$$ for $i=1,2$.

We will prove that if $r,s \in \mathbb{N}$ with $gcd(r,s)=1$ then the
ideal $\langle P_{r-1}(z,y),P_{s-1}(z,y) \rangle$ in
$\O_{z,y}$ has finite codimension. In fact $(z-y)P_{r-1}(z,y)=z^r-y^r=\prod_{i=0}^{r-1}(z-\omega^iy)$ and
$(z-y)P_{s-1}(z,y)=z^s-y^s=\prod_{i=0}^{s-1}(z-\upsilon^iy),$ with
$\omega$ (resp. $\upsilon$) $r$-th (resp. $s$-th) primitive root of
unity. Notice that $\eta$ is a $d$-th root of unity, with
$d=gcd(r,s)$ if and only if $\eta$ is a $r$-th and $s$-th root of
unity, that is, the polynomials $z^r-y^r$ and $z^s-y^s$ have $d$
common roots. In this way by hypothesis $(z-y)P_{r-1}(z,y)=(z-y)f_1$ and
$(z-y)P_{s-1}(z,y)=(z-y)f_2,$ where $f_1$ and $f_2$ don't have
common roots. Therefore $P_{r-1}$ and $P_{s-1}$ have no common
components which implies that the codimension of the ideal generated by them is finite.

Since $$d(f)=\frac{1}{2}
dim_{\Cp}\frac{\O_{\underline{x},y,z}}{\langle
	\underline{x},\varphi_{2n-1},\varphi_{2n}\rangle}=\frac{1}{2}
dim_{\Cp}\frac{\O_{y,z}}{\langle
	P_{m_1-1}(z,y),P_{m_2-1}(z,y)\rangle}$$
in which $gcd(m_1,m_2)=1$, we have that the polynomials $P_{m_i-1},i=1,2$ have no common
factor, then $d(f)$ is finite and by B\'{e}zout's Theorem we have the required formula, remembering that the delta invariant of $\phi_2$ is $\delta_{\phi_2}=\frac{(m_1-1)(m_2-1)}{2}$.

For the item $2$, the formula for $d(f)$ can be easily obtained by (\ref{FormDoubPointsQuasih}).\cqd

\subsection{$\A$-invariants under projections}

For the study of invariants of a map $f$ under some projection of it, we only assume $f \in \M_n\O(n,p)$ in which $p \geq n$ and $0$ is an isolated singularity of $f$.

Let $k\in \{1, \dots, p-1\}$ be fixed and $\pi_k:
\Cp^{p} \to \Cp^k$ be a canonical projection, and $g_k=\pi_k \circ f$ such that
$f \underset{\A}{\sim} (g_k, h),$ where
$h:(\Cp^n,0) \to (\Cp^{p-k},0)$.
Notice that $\theta(f) \simeq \theta(g_k) \oplus \theta(h)$. It is well known that
$\omega \in \Tg$ if and only if there exist
$\epsilon \in \theta(n)$ and $\eta
\in \theta(p)$ such that $\omega = df(\epsilon) + \eta(f)$.

We consider the subspaces
\begin{center}
	$T_k = \left\{ \omega_k \in \theta(g_k); \
	\omega_k= dg_k(\epsilon) + (\pi_k\circ \eta)(f), \ \epsilon \in
	\theta(n), \eta \in \theta(p) \right\}\ \mbox{and}$
	\vspace{.2cm}
	
	$T = \{\omega_{h} \in \theta(h); \
	(0,\omega_h) \in \Tg\}$\end{center}
and their codimensions are respectively defined as
$cod \ T_k = \dim_{\Cm} \displaystyle\frac{\theta(g_k)}{T_k}$ and $cod \ T = \dim_{\Cm} \displaystyle\frac{\theta(h)}{T}.$

\begin{proposition}
	\label{codT0codT}
	Let $f: (\Cp^n,0) \to (\Cp^{p},0)$ with $p\geq n$, such that $f \underset{\A}{\sim} (g_k,h)$ as above, for some $k=1, \dots, p-1$. If $cod \ T_k$ and $cod \ T$ are finite, then $f$ is $\A$-finitely determined and
	$$\Acod = cod \ T_k + cod \ T.$$
\end{proposition}
\Dem The proof is similar to the Proposition 2.1
in \cite{MRT} for divergents diagrams. If $\omega \in \theta(f)$, we denote by $\overline{\omega}$ the class of this element in the quotient
$\theta(f)/\Tg$, similarly (with an abuse of notation) for $\theta(g_k)$ and $\theta(h)$. The
result follows considering the following exact sequence:
$$0 \ \longrightarrow \ \frac{\theta(h)}{T} \ \stackrel{i^*}
{\longrightarrow} \ \frac{\theta(f)}{\Tg} \
\stackrel{\pi^*}{\longrightarrow} \ \frac{\theta(g_k)}{T_k} \
\longrightarrow \ 0,$$

\noindent where if $\overline{\omega_h} \in
\theta(h)/T$ then $i^*(\overline{\omega_h}) = \overline{(0, \omega_h)}$ and given $\omega=(\omega_1, \ldots, \omega_p)\in \theta(f)$ the projection
$\pi^*(\overline{\omega}) =
\overline{(\omega_1, \dots, \omega_k)}$.\cqd

In a natural way, we can extend the above result for multigerms.

\vspace{0.2cm}
For the germ $h:(\Cp^n,0) \to (\Cp^{p-k},0)$ we define the map $\overline{w}h: \theta(p) \to \theta(h)$ by $\overline{w}h(\eta)=(\eta_{k+1}\circ f, \dots, \eta_{p}\circ f)$, where $\eta:=(\eta_1,\ldots, \eta_p)\in \theta(p)$. Given $\omega_h \in T$ we have $(0,\omega_h) \in \Tg$, which implies that there exist
$\epsilon \in \theta(n)$ and $\eta
\in \theta(p)$ such that $dg_k(\epsilon) + (\pi_k\circ \eta)(f)=0$ and $dh(\epsilon) + \overline{w}h(\eta)=\omega_h$. Let
\begin{equation}
\label{RhoMap}
\rho:\frac{\theta(h)}{\overline{\omega}h(\theta(p))} \longrightarrow
\frac{\theta(h)}{T}
\end{equation} such that $\rho$ of the class of $\omega_h \in \theta(h)$ in $\theta(h)/\overline{\omega}h(\theta(p))$ is the class of this element in $\theta(h)/T$.

 In the next result we present a relation among $\Acod$ and invariants of some projection $g_k$ of $f$, that can be interesting if you have some information about $g_k$.

\begin{theorem}
	\label{AcodDelta}
	Consider $f: (\Cp^n,0) \to (\Cp^{p},0)$ with $p\geq 2n$ and $f \underset{\A}{\sim} (g_k,h)$. If $g_k$ is $\A$-finitely determined and $\delta_f$ is finite, then
	
	\begin{enumerate}
		\item[(i)] 	$\Acod \leq
		\A_ecod(g_k)+(p-k)\delta_{f}. $
		\item[(ii)] If $g_k$ is stable, then $\Acod =cod \ T= (p-k)\delta_{f}-dim_{\mathbb{C}} Ker(\rho)$.
	\end{enumerate}
\end{theorem}
\Dem
$(i)$ Notice that $\theta(h)/T$ is isomorphic to a submodule of $\bigoplus_{j=1}^{p-k} \O_n / f^*(\O_{p})$ and
therefore $cod \ T \leq (p-k) \delta_{f}<\infty$. On the other hand, the extended tangent space to $g_k$ is $$T{\cal A}_e(g_k) = \{\omega_k \in \theta(g_k); \omega_k= dg_k(\epsilon) +
\eta(g_k), \epsilon \in \theta(n),\eta \in \theta{(k)} \},$$ which implies $cod \ T_k \leq \A_ecod(g_k)$. By Proposition \ref{codT0codT}, $\Acod \leq
\A_ecod(g_k)+(p-k)\delta_{f}.$

\vspace{0.1cm}
$(ii)$ If $g_k$ is
stable, then $\A_ecod(g_k) = 0$ and consequently $cod \ T_k=0$. Thus $\Acod = cod \ T$.
If $\rho$ is as defined in (\ref{RhoMap}), it follows by exact sequence $$0 \ \longrightarrow \ Ker(\rho)
\ \longrightarrow \ \frac{\theta(h)}{\overline{\omega}
	h(\theta(p))} \ \stackrel{\rho}{\longrightarrow} \
\frac{\theta(h)}{T} \ \longrightarrow \ 0$$ that $cod \
T=dim_{\Cp}\frac{\theta(h)}{\overline{\omega} h(\theta(p))} -
dim_{\mathbb{C}} Ker(\rho)$. But
$\displaystyle\frac{\theta(h)}{\overline{\omega} h(\theta(p))}$
is isomorphic to $\bigoplus_{j=1}^{p-k} \frac{\O_n}{
	f^*(\O_{p})}.$\cqd

It is well known that the invariants $\delta_f$ and $d(f)$ coincide for parameterized plane curves, but this does not hold in general. In \cite[Proposition 2.9]{RodRuas} we proved that parameterized surfaces $f(x,y)=(x,x^{\lambda}y,y^{m_1},y^{m_2})$ satisfy $d(f)\leq \delta_f$ with equality if and only if $m_1=2$.

We prove in the next proposition that these two invariants coincide for any $f\in \O(n,2n)$ which has a stable projection equivalent to a fold. Moreover we obtain an estimate for $Ker(\rho)$ in Theorem \ref{AcodDelta} and consequently for the $\A_e$-codimension. For what follows, we denote the ring $\O_n$ by $\O_{\underline{x},y}$.

\begin{proposition}
	\label{DobraDelta} Let $f: (\Cp^n,0) \to (\Cp^{p},0)$ be an
	$\A$-finite map-germ for $p \geq 2n$. We suppose $f\underset{\A}{\sim}(g_n,h)$ such that
	$g_n(\underline{x},y)=(\underline{x},y^2)$. Thus
	\begin{enumerate}
		\item if $p=2n$, then $d(f)= \delta_f$;
		\item if $p\geq 2n$, then $(p-2n+1)\delta_f \leq \Acod \leq (p-n)\delta_f-1$.
	\end{enumerate}
\end{proposition}
\Dem  We can consider change of coordinates such that $f \underset{\A}{\sim}
(\underline{x},
y^2,yu_{n+1}(\underline{x},y^2),\ldots,yu_{p}(\underline{x},y^2))$
for some $u_{n+1},\ldots, u_{p} \in \Cp\{\underline{x},y\}$, as in
\cite{Mond}.

If $p=2n$, then by (\ref{d(f)G}) we have that
$\varphi_n(\underline{x},y,z)=z+y$ and
$\varphi_j(\underline{x},y,z)=\frac{zu_j(\underline{x},z^2)-yu_j(\underline{x},y^2)}{z-y}$
for all $j=n+1, \dots, 2n$. Identifying $z=-y$ we have
$\varphi_j(\underline{x},y,z)=u_j(\underline{x},y^2)$. However
$\O_{\underline{x},y} \simeq \O_{\underline{x},y^2} \oplus y
\O_{\underline{x},y^2}$ thus
$$\langle u_{n+1},\ldots, u_{2n}\rangle \O_{\underline{x},y}=\langle
u_{n+1},\ldots,u_{2n}\rangle \O_{\underline{x},y^2} \oplus \langle
u_{n+1},\ldots, u_{2n}\rangle \ y\O_{\underline{x},y^2}.$$

Hence, $d(f)=\displaystyle\frac{1}{2}dim_{\Cp}
\displaystyle\frac{\O_{\underline{x},y}}{\langle u_{n+1},\ldots,
	u_{2n}\rangle \O_{\underline{x},y}} = dim_{\Cp}
\displaystyle\frac{\O_{\underline{x},y^2}}{\langle u_{n+1},\ldots,
	u_{2n}\rangle \O_{\underline{x},y^2}}.$

\vspace{0.2cm} On the other hand, we can identify $f^*(\O_{2n})$ with
$\O_{\underline{x},y^2} \oplus \langle u_{n+1},\ldots, u_{2n}\rangle
y\O_{\underline{x},y^2}$ which implies
$$\delta_f  = \displaystyle dim_{\Cp}
\frac{\O_{\underline{x},y^2} \oplus
	y\O_{\underline{x},y^2}}{\O_{\underline{x},y^2} \oplus \langle
	u_{n+1},\ldots, u_{2n}\rangle  \ y\O_{\underline{x},y^2}} =
dim_{\Cp} \frac{\O_{\underline{x},y^2}}{\langle u_{n+1},\ldots,
	u_{2n}\rangle \O_{\underline{x},y^2}} = d(f).$$

In order to prove the item $2$, notice that it follows from Theorem \ref{AcodDelta} that $\Acod = cod \ T$, since
$g_n$ is a stable map. But
$\omega=(\omega_{n+1}, \dots, \omega_{p}) \in T$ if and only if
$(0,\omega) \in \Tg$, that is, there exist $\epsilon=(\epsilon_1,
\ldots, \epsilon_n) \in \theta(n)$ and $\eta=(\eta_1, \ldots,
\eta_{p}) \in \theta(p)$ such that $\epsilon_i=-\eta_i(f)$ for
$i=1, \ldots, n-1$ and $2y\epsilon_n=-\eta_n(f)$ which implies that for each $j$ with $n+1 \leq j \leq p$,
$$\omega_j= -\sum_{i=1}^{n-1}\left(y\frac{\partial u_j}{\partial
	x_i}({\underline
	x},y^2)\eta_i(f)\right)+\left(u_j(\underline{x},y^2)+2y^2\frac{\partial
	u_j}{\partial y}({\underline
	x},y^2)\right)\epsilon_n(\underline{x},y)+\eta_j(f)$$ satisfies $\gamma_j:= u_j({\underline
	x},y^2)\epsilon_n({\underline x},y)-y \frac{\partial u_j}{\partial
	y}({\underline x},y^2)\eta_n(f)+\eta_j(f)\in f^*(\O_{p})$. In fact,
$y\frac{\partial u_j}{\partial y}({\underline x},y^2)$ is even in
$y$ and by isomorphism $\O_{\underline{x},y} \simeq
\O_{\underline{x},y^2} \oplus y \O_{\underline{x},y^2}$ we can write
$\epsilon_n(\underline{x},y)=\epsilon_{n_1}(\underline{x},y^2)+y
\epsilon_{n_2}(\underline{x},y^2)$. Thus $$u_j({\underline
	x},y^2)\epsilon_n({\underline x},y)=u_j({\underline
	x},y^2)\epsilon_{n_1}(\underline{x},y^2)+yu_j({\underline
	x},y^2)\epsilon_{n_2}(\underline{x},y^2)\in f^*(\O_{p})$$ because
$u_j,\epsilon_{n_1},\epsilon_{n_2}$ are even in $y$ and
$yu_j({\underline x},y^2)\in f^*(\O_{p})$. It follows from (\ref{RhoMap}) and the proof of Theorem \ref{AcodDelta} that $\omega \in
\frac{T}{\overline{\omega}h(\theta(p))}$ if and only if $\omega \in M$, where $$M:=\frac{f^*(\O_{p})\{v_1, \dots, v_{n-1}\}
	+\overline{\omega}h(\theta(p))}{\overline{\omega}h(\theta(p))}$$  is a
$f^*(\O_{p})$-module finitely generated by
$$v_i=\left(\begin{array}{c} y\frac{\partial u_{n+1}}{\partial
	x_i}\\
\vdots \\
y\frac{\partial u_{p}}{\partial x_i}
\end{array}\right),$$ with $i=1, \ldots, n-1$. By abuse of notation we denote $v_i$ its class in
$\overline{\omega}h(\theta(p))$.

Notice that $v_1, \dots, v_{n-1}$ are not all identically zero in
$M$. In fact, since $p\geq 2n$ and $f$ is $\A$ (and therefore ${\cal
	L}$) finitely determined, for each $i=1, \ldots, n-1$ there exists
$k_i \in \mathbb{N}$ such that $x_i^{\gamma_i}y \in f^*(\O_{p})$ for
all $\gamma_i \geq k_i$. This implies that writing
$u_{n+i}(\underline{x},y^2)=\phi_{n+i}(\underline{x})+y^2R_{n+i}(\underline{x},y)$ there are at least $n-1$ maps $\phi_j$'s not identically
zero, that we can suppose they are $\phi_{n+1}, \ldots , \phi_{2n-1}$. For each $i=1, \ldots, n-1$ there exists $l_i$ such that
$\phi_{l_i}(\underline{x})=a_{l_i}x_i^{l_i}+\ldots$ with
$ord_{x_i}\phi_{l_i}(\underline{x}) \leq
ord_{x_i}\phi_{j}(\underline{x})$ for all $j$ with $l_i,j  \in
\{n+1, \ldots, 2n-1\}$. Thus $$ord_{x_i}\frac{y\partial
	u_{l_i}(x,y^2)}{\partial x_i} < l_i$$ and, therefore, the
$l_i$-coordinate of $v_i$ not belongs to $f^*(\O_{p})$, by minimality
of the order of $\phi_{l_i}$ with respect to $x_i$. Consequently
$v_i \not\in \overline{\omega}h(\theta(p))$ and we conclude that
$dim_{\Cp}\ Ker(\rho) \geq 1$.

If $v_1, \ldots, v_{n-1}$ are non zero in $M$, this means that each
non-zero coordinate of $v_i$ is an element of
$\frac{\O_n}{f^*(\O_{p})}\simeq \Cp\{\lambda_1, \ldots,
\lambda_{\delta_f}\}$. Since $f$ is ${\cal L}$-finite, for each $i=1, \ldots, n-1$ there exist (finite)
$m_{i1},\ldots, m_{ir_i}\in f^*(\O_{p})$ such that $m_{ij}v_i
\not\in \overline{\omega}h(\theta(p))$, with $1 \leq j \leq r_i$.
But, in each non zero component of $v_i, m_{i1}v_i,\ldots,
m_{ir_i}v_i$ we have at most $\delta_f$ linearly independent
elements in $\frac{\O_n}{f^*(\O_{p})}$. According to these
arguments we have
$dim_{\Cp} M \leq (n-1)\delta_f$. Therefore $1\leq dim_{\Cp}\ Ker(\rho)=dim_{\Cp}\
\frac{T}{\overline{\omega}h(\theta(p))}\leq (n-1)\delta_f$. The
result follows from Theorem \ref{AcodDelta}. \cqd

\begin{corollary}
	\label{CorDobraDelta}
	If $v_1, \ldots, v_{n-1}$ (as defined in the above proposition) are
	$f^*(\O_{p})$-linearly independent, then $\Acod \leq (p-n) \delta_f
	-n+1$.
\end{corollary}
\Dem If $v_1, \ldots, v_{n-1}\in M$ are $f^*(\O_{p})$-linearly
independent we have $dim_{\Cp}\ Ker(\rho) \geq n-1$.\cqd

\begin{example} In this example we present $\A$-simple map-germs from the
	classification obtained in \cite{Klotz}. Notice that the estimates for the $\A_e$-codimension in Proposition \ref{DobraDelta} are sharp. In fact, the germs $f(x,y)=(x,xy,y^2, y^{2k+1})$ and $g(x,y)=
	(x,y^2,y^3,x^ky)$ satisfy $\delta_f=k=\delta_g$. Thus $\Acod = k=\delta_f$ and $\A_e
	cod(g)=2k-1=2\delta_g -1$.
	
	Now, the map $f(x_1,x_2,y)=(x_1,x_2,y^2,x_1x_2y,(x_1^2-x_2^2)y,y^3)$ illustrates the computation of $Ker(\rho)$. It is possible to prove that $\frac{\O_3}{f^*(\O_{6})}\simeq \Cp\{y,x_1y,x_2y,x_1^2y\}$, thus
	$\delta_f= 4$. In this case, $$Ker(\rho)\simeq \displaystyle\frac{\Cp\{v_1,v_2,x_1v_1,x_1v_2\}
		+\overline{\omega}h(\theta(6))}{\overline{\omega}h(\theta(6))}$$ is a $f^*(\O_{6})$-module
	generated by $\{v_1,v_2 \}$, where $$v_1= \left(\begin{array}{c}
	x_1y \\
	-2x_2y\\
	0
	\end{array}\right) \ \ \mbox{and} \ \ v_2=\left(\begin{array}{c}
	x_2y \\
	2x_1y\\
	0
	\end{array}\right).$$ By Theorem \ref{AcodDelta}, we conclude that $\Acod = 3 \delta_f - dim_{\Cp}Ker(\rho)=12-4
	=8 < 3 \delta_f - n+1.$
\end{example}

\vspace{0.5cm}

\noindent Rodrigues Hernandes, M. E. \\
{\it merhernandes@uem.br} \\
Universidade Estadual de Maring\'{a}, DMA, Av. Colombo 5790.
Maring\'{a}-PR 87020-900,  Brazil.

\vspace{0.4cm}
\noindent Ruas, M. A. S. \\
{\it maasruas@icmc.usp.br} \\
Universidade de S\~{a}o Paulo, ICMC, Caixa Postal 668, S\~{a}o Carlos-SP, 13560-970, Brazil.

\begin{thebibliography}{99}
\bibitem{BRS} Benedini Riul, P., Oset Sinha, R. and Ruas, M. A. S., \emph{The geometry of corank $1$ surfaces in $\mathbb{R}^4$}.
Quart. J. Math. \textbf{70} (2019), 767--795.

\bibitem{BMN} Birbrair, L., Mendes, R. and Nuño-Ballesteros, J. J., \emph{Metrically un-knotted corank 1 singularities of surfaces in $\mathbb{R}^4$}. J. Geom. Analysis \textbf{28} (2018), 3708--3717.

\bibitem{CMW} Cooper, T., Mond, D. and Wik Atique, R., \emph{Vanishing topology of codimension $1$ multi-germs over $\mathbb{R}$ or $\mathbb{C}$}.
Comp. Math. \textbf{131} (2002), 121--160.

\bibitem{ES} Eisenbud, D. and Sturmfels, B., \emph{Binomial ideals.} Duke Math. J. \textbf{84} (1996), no. 1, 1--45.

\bibitem{Gaffney} Gaffney, T., \emph{${\cal L}^0$-equivalence of
maps.} Math. Proc. Camb. Phil. Soc. \textbf{128} (2000), 479--496.

\bibitem{GolubitskyGuil} Golubitsky, M. and Guillemin, V., \emph{Stable Mappings and their Singularities.} Graduate Texts in Mathematics 14, NY.
Springer-Verlag (1973).

\bibitem{Greuel} Greuel, G.M., \emph{Equisingular and equinormalizable deformations of isolated non-normal singularities.} Methods Appl. Anal. \textbf{24} (2017), no. 2, 215--276.

\bibitem{Hefez} Hefez, A., \emph{Irreducible Plane Curve
Singularities.} Real and Complex Singularities. Eds D. Mond and M. J.
Saia. Lecture Notes in Pure and Appl. Math. 232, NY. Marcel
Dekker (2003).

\bibitem{Ishikawa} Ishikawa, G. and Janeczko, S., \emph{Symplectic invariants of parametric singularities.}
Advances Geom. Analysis, ALM, \textbf{21} (2011), 259--280.

\bibitem{Klotz} Klotz, C., Pop, O. and Rieger, J. H., \emph{Real double-points of deformations of
$\A$-simple map-germs from $\mathbb{R}^n$ to $\mathbb{R}^{2n}$.}
Math. Proc. Camb. Phil. Soc. \textbf{142} (2007), 341--363.

\bibitem{MRT} Mancini, S., Ruas, M. A. S. and Teixeira, M. A., \emph{On divergent diagrams of finite codimension.}
Portugaliae Math.  \textbf{59} (2002), no. 2, 179--194.

\bibitem{Mond} Mond, D., \emph{On the classification of germs of maps from $\mathbb{R}^2$ to $\mathbb{R}^3$}.
Proc. London Math. Soc., \textbf{50} (1985), no. 3, 333--369.

\bibitem{ORW} Oset Sinha, R., Ruas, M. A. S. and Wik Atique, R., \emph{Classifying codimension two multigerms}. Math. Z. \textbf{278} (2014), 547--573.

\bibitem{P} Peñafort-Sanchis, G., \emph{Reflection maps}.
Math. Annalen \textbf{378} (2020), 559--598.

\bibitem{RodRuas} Rodrigues Hernandes, M. E. and Ruas, M. A. S., \emph{Parametrized monomial surfaces in $4$-space.} Quart. J. Math. \textbf{70} (2019), 473--485.

\bibitem{Wall} Wall, C. T. C.,  \emph{Finite determinacy of smooth map-germs.} Bull. London Math. Soc.
\textbf{13} (1981), no. 6, 481--539.

\bibitem{Whitney} Whitney, H.,  \emph{On singularities of mappings of Euclidean spaces. I. Mappings of the plane into the plane}. Annals Math. 
\textbf{62} (1955), no. 3, 374--410.

\end{thebibliography}
\end{document}